\documentclass[journal]{IEEEtran}

\usepackage{fancyhdr}

\usepackage{setspace}
\usepackage{amsmath, amssymb, amsfonts}
\usepackage{enumerate}
\usepackage{graphicx}

\usepackage{verbatim}
\usepackage{url}

\usepackage{graphicx}
\usepackage{color}
\usepackage{amsthm}
\usepackage{amsmath}
\usepackage{multicol}

\newcommand{\Int}{\operatorname{int}}

\newcommand{\diag}{\operatorname{diag}}
\newcommand{\sign}{\operatorname{sgn}}

\newcommand*\diff{\mathop{}\!\mathrm{d}}

\newcommand{\st}{\, | \,}

\newtheorem{Theorem}{Theorem}

\newtheorem{Proposition}{Proposition}
\newtheorem{Lemma}{Lemma}
\newtheorem{Corollary}{Corollary}
\newtheorem{Definition}{Definition}
\newtheorem{Example}{Example}

\newtheorem{Remark}{Remark}

\newcommand {\R}{\mathbb R}

\renewcommand {\S}{\mathbb S}

\newcommand {\D}{\mathbb D}

\newcommand{\be}{\begin{equation}}
\newcommand{\ee}{\end{equation}}
\newcommand{\sgn}{\operatorname{{\mathrm sgn}}}
\newcommand{\updt}[1]{{\color{black}#1}}

 \IEEEoverridecommandlockouts
\begin{document}
%%%%%%%%%%%%%%%%%%%%%%%%%
\title{Diagonal Stability of Discrete-time $k$-Positive linear Systems with Applications to Nonlinear Systems  \thanks{
		Research  supported in part 
		by  research grants from    
		the Israel Science Foundation and the 
		US-Israel Binational Science Foundation.  
}}

%%%%

\author{Chengshuai Wu and Michael Margaliot\thanks{The authors
 are with  the School of EE-Systems, Tel Aviv University, Israel 69978. Corresponding author: Michael Margaliot,
Email: michaelm@tauex.tau.ac.il}}

\maketitle 

 %%\onehalfspace

\begin{abstract}
%%%%%%%%%%%%%%%%%%%%%%%%%%
A linear dynamical system is called 
 $k$-positive if its dynamics maps the set of vectors with up to~$k-1$ sign variations to itself. For~$k=1$, this reduces to the important class  of positive linear systems. 
Since stable positive linear time-invariant (LTI) systems always admit a \emph{diagonal} quadratic Lyapunov function, i.e. they are diagonally stable, we may expect
that this holds also for stable~$k$-positive systems. We show that, in general, this is not the case both in the continuous-time~(CT) and discrete-time~(DT) case.
We then focus on  DT $k$-positive linear systems and introduce the new notion of
 \emph{DT $k$-diagonal stability}. It is shown that this 
 is a necessary condition for   standard~DT diagonal stability.
We demonstrate an application of this new notion 
to the   analysis of  a class of DT nonlinear systems.
%%%%%%%%%%%%%%%%
\end{abstract}


%%%%%%%%%%%%%%%%%%%%%%%%%%%%%%%%%
\section{Introduction}
%%%%%%%%%%%%%%%%%%%%%%%%%%%%%%%%%
 Lyapunov functions are a powerful tool for stability  analysis
and control synthesis. For linear time-invariant~(LTI) systems, stability is equivalent to the existence of a quadratic Lyapunov function, i.e.~$V(x)=x^T Q x$, with~$Q$ positive-definite, 
that can be obtained constructively based on the eigenvectors of an associated Hamiltonian matrix~\cite{linear_st}. An~LTI is called \emph{diagonally  stable} if it is possible to find a diagonal Lyapunov function~(DLF), i.e.~$V(x)=x^T D x$, with~$D$ 
 positive-definite and diagonal.
%%%%%%%%%%%%%%%%%%%%%%%%%%%%%%%%%%

Diagonal stability of LTIs has attracted considerable attention in the systems and control community (see e.g. the monograph~\cite{kaszkurewicz2012matrix}).
Due to its simplicity, diagonal stability can facilitate control synthesis, and it plays an important role in many fields including mathematical 
economics~\cite{barker1978positive}, ecology~\cite{cross1978three}, numerical analysis~\cite{kraaijevanger1991unique}, biochemistry~\cite{arcak2006diagonal}, and networked systems~\cite{rantzer2015scalable}. 

The existence of a DLF has important implications to certain \emph{nonlinear} systems associated with the~LTI~\cite{kaszkurewicz1984note, kaszkurewicz1993robust}.
This is true for both continuous-time~(CT) and discrete-time~(DT) nonlinear systems. 
We now briefly explain this. For~$P\in \R^{n\times n}$, we write~$P\succ 0$ [$P\prec 0$]
to denote that~$P$ is symmetric and positive-definite [negative-definite].
%%%%%%%%%%%%%%%%%%%%%%%%%%%%%%%%%
Consider the~CT nonlinear system:
\begin{align}  
\dot{x}(t) =& A f(x(t)),   \label{eq:nonli_ct}
\end{align}
where $A \in \mathbb{R}^{n \times n}$,
  $f(x)  = \begin{bmatrix} f_1(x_1)& \cdots& f_n(x_n) \end{bmatrix}^T$, $f_i$ is continuous and~$f_i(z) z> 0$  for all $z \neq 0$ (so~$f_i(0)=0$). 
	\updt{Such a dynamics is called a Persidskii system
	(see, e.g.~\cite{prese2019,prese2021} and the references therein).}
	%%%
	Suppose that~$A$  satisfies the Lyapunov inequality $ DA + A^T D \prec 0$ with a \emph{diagonal}
	matrix~$D \succ 0$.
	Let
	\[
	V(z):= 2\sum_{i=1}^{n} d_i \int_{0}^{z_i } f_i(\tau) \diff \tau,
	\]
	where~$d_i$ is the $i$th diagonal entry of~$D$.
	Then the derivative of~$V(x(t))$ along solutions of~\eqref{eq:nonli_ct} is
	\[
	\dot V(x(t))=  f^T(x(t))( DA + A^T D) f(x(t)),	\]
	so~$\dot V(x(t))<0$ whenever~$x(t)\not =0$.
If~$\int_0^{x_i} f_i(\tau) \diff \tau  \to \infty$ as  $\vert x_i \vert \to \infty$, $i = 1, \dots, n$,
then we can conclude that the nonlinear system~\eqref{eq:nonli_ct}  is
globally asymptotically stable~(GAS).

Note that~\eqref{eq:nonli_ct} can also be interpreted as a networked system. Indeed, 
assume that~$A$ is nonsingular and let~$y := A^{-1} x$. Then~\eqref{eq:nonli_ct}  becomes 
\begin{equation} \label{eq:net_y}
\dot{y}_i(t) = f_i \left(\sum_{j=1}^n a_{ij} y_j(t) \right),~ i = 1, \dots, n,
\end{equation}
where $a_{ij}$ is the $(i,j)$-th entry of~$A$. 
This can be viewed as a networked
 system with the weighted adjacency matrix~$A$. In this case, diagonal stability of the LTI 
implies GAS of an associated nonlinear networked system.
This idea was used in~\cite{arcak2006diagonal} to show that diagonal stability of a cyclic LTI implies the stability 
of a cyclically interconnected network of output strictly passive systems~\cite{van2017L2}.

A similar construction holds for the DT nonlinear system:
\begin{align}  
x(j+1) =& A \phi(x(j)),   \label{eq:nonli_dt}  
 \end{align}
where   $\phi(x): = \begin{bmatrix} \phi_1(x_1)& \cdots& \phi_n(x_n) \end{bmatrix}^T$, 
with~$\phi_i(z)$ continuous and~$
0 < \vert \phi_i(z) \vert \leq \vert z \vert   
$ 
   for all~$z \neq 0$ (so~$\phi(0)=0$). If~$A$ in~\eqref{eq:nonli_dt} satisfies the Stein inequality $A^T D A \prec D$, with a \emph{diagonal} matrix~$D\succ 0$, then 
	$
			V(z) : = z^T   D z	
	$
is a Lyapunov function for the nonlinear system~\eqref{eq:nonli_dt}. Similar to the CT case, the DT nonlinear system \eqref{eq:nonli_dt} can also be interpreted 
 as a networked system   using a suitable
 change of coordinates.

Stable LTIs always admit a quadratic Lyapunov function, but not necessarily  a~DLF~\cite{barker1978positive}.   
It is well-known however that   stable \emph{positive} LTIs do admit a DLF (see, e.g.,~\cite{rantzer2015scalable}).

Recently, the notion of positive linear systems was generalized to $k$-positive linear systems.
\updt{For the theory and applications of such systems, see~\cite{Eyal_k_posi, alseidi2019discrete} and also~\cite{rami_osci,is_my_system,grussler_tf,grusseler_btrunctaion}}.
%%%%%
For~$k=1$, this reduces to   positive linear systems. This naturally raises the question of whether stable $k$-positive LTIs also admit a~DLF. Here, we show that the answer is in general no.

We then focus on the DT case. We show that $k$-positive DT LTI always satisfy
 a   property \updt{that} we call \emph{DT $k$-diagonal stability}.
 It is showed that DT $k$-diagonal stability is a necessary condition
for $DT$ diagonal stability. 
We then describe an 
  application to a class of DT nonlinear systems in a form similar to~\eqref{eq:nonli_dt}. By using   wedge products and their   
	geometric interpretation, we show that the asymptotic behavior
	of these   systems can be analyzed using $k$-positivity and DT $k$-diagonal stability.  
	These result generalize the construction described above when~$k=1$. 

The remainder of this note is organized as follows. The next section briefly reviews some basic definitions and known results from the theory of
 diagonal stability, positive LTIs, compound matrices, and~$k$-positive systems.
Section \ref{sec:doesnot} shows that in general stable~$k$-positive systems, with~$k>1$, are \emph{not} diagonally stable.  
Section \ref{sec:k-diag} introduces the notion of DT $k$-diagonal stability, and explains its relation to the standard DT diagonal stability. 
An application to DT \emph{nonlinear} systems is described in Section~\ref{sec:non}.
 
We use standard notation. 
A matrix   $X \in \R^{n \times m}$
is called non-negative [positive], denoted~$X \geq 0$ [$X\gg 0$], if all its entries are non-negative [positive]. The determinant   of~$A\in \mathbb{R}^{n \times n}$ is denoted
	by~$\det(A)$. 
%%%%%%%%%%%%%%
% $A  $ is called
% \emph{Metzler} if all its off-diagonal entries are non-negative. 
%%%%%%%%%%%%%%
 The eigenvalues of $A $ are denoted 
by~$\lambda_i(A)$, $i = 1, \dots, n$, ordered such that
\be \label{lam_i}
|\lambda_1(A)|\geq|\lambda_2(A)|\geq\dots\geq |\lambda_n(A)|.
\ee
The \emph{spectral radius} of~$A$ is~$\rho(A):=|\lambda_1(A)|$. 
For two integers~$i\leq j$, we let~$[i,j]:=\{i,i+1,\dots,j\}$. 
The non-negative orthant in~$\R^n$ is
~$
\R^n_+:=\{x\in\R^n \st x_i \geq 0,\; i=1,\dots,n\}
$.

%%%%%%%%%%%%%%%%%%%%%%%%%%%%%%%%%
\section{Preliminaries} \label{sec:pre}
%%%%%%%%%%%%%%%%%%%%%%%%%%%%%%%
In this section, we review several known topics that are needed later on.

\subsection{Diagonal stability of positive DT  LTIs}
%%%%%%%%%%%%%%%%%%%%%%%%%%%%%%%%%%%%%%%%%
If $A \in \mathbb{R}^{n \times n}$ is non-negative, 
then~$x(j+1)=Ax(j)$
is called a positive  DT  LTI.
%%%%%%
% If $A \in \mathbb{R}^{n \times n}$ is Metzler, then 
% $\dot{x}(t) =  Ax(t) $
%is called  a positive continuous-time LTI (CT-LTI). 
%%%%%
The dynamics  of
positive DT  LTIs leaves the proper cone~$\R^n_+ $ invariant~\cite{hlsmith}. 
The following result  shows
 that stable positive DT  LTIs are diagonally stable. 
\updt{Let~$\D^{n\times n}$ denote the set of~$n\times n$ positive
diagonal   matrices.}
%%%%%%%%%%%%%%%%%%%%%%%%%%%%%%%%%%
\begin{Lemma} [{see e.g. \cite[Prop. 2]{rantzer2015scalable}}] \label{le_diag}
 If $A \in \mathbb{R}^{n \times n}$ with~$A\geq 0$ then
 the following statements are equivalent:
	\begin{enumerate}[(a)]
	\item The matrix $A$ is Schur, i.e., $\rho(A) < 1$;
	\item There exists $\xi \in\R^n$ with~$\xi \gg 0$ such that  $A\xi \ll \xi$;
	\item There exists  $z\in\R^n$ with~$z \gg 0 $ such that    $A^T z \ll z$;
	\item There exists  $D \in \D^{n \times n}$ such that~$A^T D A \prec D$;\label{item:dsddt}
	\item The matrix $(I - A)$ is nonsingular and $(I - A)^{-1} \geq 0 $.
	\end{enumerate} 
\end{Lemma}

\begin{Remark}[see e.g.~\updt{\cite{rantzer2015scalable}}] \label{re_proce}
Let~$A  \in  \mathbb{R}^{n \times n}$ be non-negative and Schur. 
 Pick~$x,y \in \mathbb{R}^n$ with $x,y\gg 0$. Then $\xi := (I - A)^{-1}x  $, $z := (I - A^T)^{-1}y  $, and $D := \operatorname{diag}(\frac{z_1}{\xi_1}, \dots,\frac{z_n}{\xi_n})$ satisfy  conditions (b), (c), and (d) in Lemma \ref{le_diag}, respectively. This provides a constructive procedure to obtain a~DLF for positive~DT  LTIs. Note that if~$A \in \mathbb{R}^{n \times n}$ is Schur and~$A\leq 0$, then~$(-A)$ is Schur and non-negative. In this case, Lemma~\ref{le_diag} also guarantees    the existence
of  a~$D \in \D^{n \times n}$ such that~$A^T D A \prec D$. 
\end{Remark}

%%%%%%%%%%%%%%%%%%%%%%%%%%%%%%%%%%%%%
\subsection{$k$-positive systems}
%%%%%%%%%%%%%%%%%%%%%%%
 \updt{We recall two definitions for the number of sign variations in a vector. 
Define~$s^-,s^+:\R^n \to \{0,1\dots,n-1\}$ as follows. First,~$s^-(0)=0$. Second, for~$x\not =0 $,
 $s^-(x)$ is the number of sign variations   
in~$x$ after deleting all its zero entries. 
Let~$s^+(x)$ denote the maximal possible number of sign variations
in~$x$ after each zero entry is replaced by either~$1$ or~$-1$. 
For example, for~$n=4$ and~$x=\begin{bmatrix} 1.3& 0 & 0  &  -\pi  \end{bmatrix}^T$, we have~$s^-(x)=1$ and~$s^+(x)=3$. 
Obviously,
\[
%%%
0\leq s^-(x) \leq s^+(x)\leq n-1 \text{ for all } x\in\R^n.
\]
For any~$k\in [1,n]:= \{1,\dots,n\}  $, define the sets:
\be\begin{aligned}\label{eq:defpk0}
	P^k_- &:=\{ x\in\R^n : s^-(x)\leq k-1 \}, \\
	P_+^k &:=\{ x\in\R^n : s^+(x)\leq k-1 \}.
\end{aligned}\ee
%%%%%%%%%%%%
For example,~$P^1_-=\R^n_+ \cup (-\R^n_+)$.

A linear dynamical system is called~\emph{$k$-positive} if its flow maps~$P^k_-$ to~$P^k_-$, and \emph{strongly~$k$-positive} if its flow maps~$P^k_-\setminus\{0\}$ to~$P_+^k$ \cite{Eyal_k_posi, alseidi2019discrete}.
}
For example, 
the dynamics of positive LTIs maps the non-negative orthant 
$
\R^n_+ 
$ to itself (and also $-\R^n_+$ to itself), so they are~$1$-positive systems.  

\updt{Multiplying a vector by a non-zero scalar does not change the number of sign variations in the vector. This implies that~$P^k_-,P^k_+$ are cones. 
However, they are not convex cones. For example, 
the vectors~$x:=\begin{bmatrix} 4&2&4 \end{bmatrix}^T$
and~$y:=\begin{bmatrix} -2&-4&-2 \end{bmatrix}^T$
satisfy~$x,y \in P^1_-$ and~$x,y \in P^1_+$, but~$z:=(x+y)/2=\begin{bmatrix} 1&-1&1 \end{bmatrix}^T$
satisfies~$z\not \in P^1_-$
and~$z \not \in P^1_+$.
}

The analysis of $k$-positive  systems is based on compound matrices.

%%%%%%%%%%%%%%%%%
\subsection{Multiplicative compound matrices} \label{subsec:com}
%%%%%%%%%%%%%%%%%%%%%%%%%%%%%%%%%%%%%%%%%
 For an integer~$n\geq 1$ and~$k\in [1,n]  $,
let~$Q_{k,n}$ denote the ordered set of all strictly increasing sequences of~$k$ integers
 chosen from~$[1,n]$. 
We denote the~$r:=\binom{n}{k}$  elements of~$Q_{k,n}$ by
$\kappa_1,\dots,\kappa_r$,
with the~$\kappa_i$s ordered	 lexicographically.
For example,~$Q_{2,3}=\{  \kappa_1,\kappa_2,\kappa_3    \}$,
with~$\kappa_1=\{1,2\}$, $\kappa_2=\{1,3\}$,
and~$\kappa_3=\{2,3\}$.  
%%%%%%%%%%%%%%%%%%%

Given $A \in \mathbb{R}^{n \times n}$
and~$\kappa_i , \kappa_j \in  Q_{k,n} $,
 let~$A[\kappa_i \vert \kappa_j]\in \mathbb{R}^{k \times k}$ denote the submatrix of~$A$ consisting of the rows [columns] indexed by~$\kappa_i$ [$\kappa_j$].  
Let
$
A(\kappa_i \vert \kappa_j) := \det(A[\kappa_i \vert \kappa_j]),
  $
i.e., the $k$-minor of $A$ determined by the rows [columns] in~$ \kappa_i$ [$ \kappa_j$]. 
%%%%%%%%%%%%%%%%%%%%%%%%

The \emph{$k$th multiplicative compound~(MC)}
of~$A$ is the matrix~$A^{(k)} \in \mathbb{R}^{r \times r}$, whose entries, written in lexicographic order, are $A(\kappa_i \vert \kappa_j)$, see e.g.~\cite{margaliot2019revisiting,Eyal_k_posi} for more detailed explanations and examples. Note that this implies that~$A^{(1)}=A$ and~$A^{(n)}=\det(A)$. 
The~MC satisfies the following properties (see, e.g.,~\cite{muldo1990}). 
%%%%%%%%%%%%%%%%%%%%%%%%%%%%%%%%%%%%%%%%%%%%%%%%%%%%%
\begin{Lemma} \label{le_multi}
Let~$A,B\in\R^{n\times n}$ and pick~$k\in [1,n]$. 
Then
\begin{enumerate}[(a)]
	\item   $(AB)^{(k)} = A^{(k)}B^{(k)}$;\label{item:i}
	\item \label{item:ii} if~$A$ is nonsingular then $(A^{-1})^{(k)} = (A^{(k)})^{-1}$;
	\item \label{item:iii} $(A^T)^{(k)} = (A^{(k)})^T$;
	\item  \label{item:iv} if~$A^{\frac12} $ exists 
	then $(A^{\frac12})^{(k)} = (A^{(k)})^\frac12$;
	\item \label{item:eig} the product of every $k$ eigenvalues of~$A$ is an eigenvalue of~$A^{(k)}$;
	\item \label{item:v} if $A$ is Schur, then $A^{(k)}$ is Schur; 
	\item \label{item:vi} if $A$ is a diagonal matrix, then $A^{(k)}$ is a diagonal matrix.
	\item \label{item:vii} if $A \succ 0$, then~$A^{(k)} \succ 0$.
\end{enumerate}
\end{Lemma}

Note that Property~\eqref{item:i} justifies the term  multiplicative compound. For~$k=n$, this property 
becomes the familiar formula~$\det(AB)=\det(A)\det(B)$.
 
 %%%%%%%%%%%%%%%%%%%%%%%%%%%%%%%%%%%%%%%%%%%%%%%%%%%%

\subsection{Necessary and sufficient conditions for $k$-positivity}
%%%%%%%%%%%%%%%%%%%%%%%%%%%%%%%%%%%%%%%
\updt{A matrix~$A\in \R^{n\times m}$
is called  \emph{sign-regular of order~$k$}, denoted~$SR_k$, if either $A^{(k)} \leq 0$ or $A^{(k)} \geq 0$. It   
is called  \emph{strictly sign-regular of order~$k$}, denoted~$SSR_k$, if either~$A^{(k)} \ll 0$ or~$A^{(k)} \gg 0$.}
In other words, all minors of order~$k$ 
of~$A$ have the same [strict] sign.\footnote{We note that the terminology in this field is not uniform and some authors refer to such matrices as \emph{sign-consistent of order $k$}.}
To refer to the common sign of the entries of $A^{(k)}$, we use the \emph{signature}~$\epsilon_k \in \{-1,1\}$. That is, if~$A^{(k)}$ is $SSR_k$ [$SR_k$]
with  signature~$\epsilon_k=1$, then all the~$k$-minors of~$A$ are positive [non-negative].
%%%%%%%%%%%%%%%%%%%%%%%%%%%%%

The next result provides a necessary and sufficient condition for a nonsingular matrix to map~$P^k_-$ to itself.
%%%%%
\begin{Proposition}[\cite{CTPDS}] \label{prop:sr}
	Let $T \in \mathbb{R}^{n \times n}$ be a nonsingular matrix and pick $k \in [1,n]$. Then 
	\begin{enumerate}[(a)]
		\item   $T P^k_- \subseteq  P^k_-$ if and only if $T$ is $SR_k$;
		\item   $T ( P^k_-\setminus\{0\}) \subseteq P^k_+$ if and only if $T$ is $SSR_k$. \label{enu:thm1b}
	\end{enumerate} 
\end{Proposition}

For example, for~$k=1$ this implies that~$T (\R^n_+ \cup (-\R^n_+)) \subseteq (\R^n_+ \cup (-\R^n_+))$ if and only if~(iff)
    the entries of~$T$ are all non-negative or all non-positive,
and that~$T (\R^n_+ \cup (-\R^n_+)) \subseteq \Int(\R^n_+ \cup (-\R^n_+))$
iff   the entries of~$T$ are all  positive or all negative. 

\begin{Remark}
The assumption that~$T$ is nonsingular is not restrictive in our setting.
Indeed, if~$x(j +1 ) = Ax(j)$, with~$A$ singular, then the dynamics can be reduced to a lower-dimensional DT  LTI   with a nonsingular matrix.   
\end{Remark}

%%%%%%%%%%%%%%%%%
 The next result gives a necessary and sufficient condition for 
a DT  LTI to be~$k$-positive. 
%%%%%%%%%%%%%%%%%%%%%
\begin{Proposition}[\updt{\cite[Thm. 1]{alseidi2019discrete}}] \label{prop:dt-lit}
	Let $A \in \mathbb{R}^{n \times n}$ be   nonsingular  and pick~$k \in [1,n]$. 
	The DT  LTI
	\be \label{dt_lti}
	x(j+ 1) = Ax(j)
	\ee
	is  $k$-positive iff $A $ is~$SR_k$, and strongly $k$-positive iff $A $ is~$SSR_k$.
\end{Proposition}

\subsection{Wedge products} \label{subsec:wed}
%%%%%%%%%%%%%%%%%%%%%%%%%%%%%%%%%%%%%%%%%%%%%%%%%%%%%%%%%%%%%%%%%%%%%%
Fix an integer~$n\geq 1$ and~$k\in[1,n]$. The wedge product of the~$k$ vectors~$a^1,\dots,a^k \in \R^n$ is
defined as
\be\label{eq:wdfe}
a^1\wedge\dots\wedge a^k:=\begin{bmatrix} a^1 & \dots &a^k \end{bmatrix}^{(k)}.
\ee
We also use the notation~$\wedge_{i=1}^k a^i := a^1\wedge\dots\wedge a^k$. Note that the right-hand side of \eqref{eq:wdfe} has dimensions~$\binom{n}{k} \times \binom{k}{k}$, that is, it is 
a column vector of dimension~$\binom{n}{k}$. 
In the special case~$k=n$, Eq.~\eqref{eq:wdfe} yields
\[
\wedge_{i=1}^n a^i =\begin{bmatrix} a^1 & \dots &a^n \end{bmatrix}^{(n)}=\det(\begin{bmatrix} a^1 & \dots &a^n \end{bmatrix}).
\]

\begin{comment}
Let~$A\in\R^{n\times n}$. Then  
\begin{align*}
   \wedge_{i=1}^k A	a^i	&	=	\begin{bmatrix} Aa^1 & \dots &Aa^k \end{bmatrix}^{(k)}		\\
    &=(A\begin{bmatrix}  a^1 & \dots & a^k \end{bmatrix})^{(k)}	\\
    &=A^{(k)}\begin{bmatrix}  a^1 & \dots & a^k \end{bmatrix} ^{(k)}\\
	&=A^{(k)} \wedge_{i=1}^k a^i.
\end{align*}
\end{comment}

The wedge product has an important geometric meaning. The value~$|\wedge_{i=1}^k a^i|$  
is the \emph{$k$-content} \cite{muldo1990} of the 
parallelotope whose edges are the given vectors. For $k =2$ and $k=3$, the~$k$-content reduces to the standard notion of area and volume.
For example, consider the case~$n=3$ and~$k=2$. Pick~$a,b \in \R^3$. Then
\begin{align*}
								a\wedge  b	&	=	\begin{bmatrix} a_1 & b_1\\a_2 & b_2\\a_3 & b_3  \end{bmatrix}^{(2)}\\		
								 =&\begin{bmatrix} a_1 b_2 -b_1 a_2 & a_1 b_3-b_1 a_3 & a_2 b_3-b_2 a_3 \end{bmatrix} ^T	.
%%%%
\end{align*}
 The entries here  are the same as those in the   cross product~$a \times b$, up to a minus sign.
Thus,~$|a \wedge b|=|a\times b|$, and when~$|\cdot|$ is the Euclidean norm this is the area  
 of the parallelogram having~$ a $ and~$ b$ as sides.

\subsection{Necessary conditions for diagonal stability} \label{subsec:nec_ds}
%%%%%%%%%%%%%%%%%%%%%%%%%%%%%%%%%%%
Recall    that~$A(\kappa_i\vert \kappa_j)$ is 
 a \emph{principal minor} of~$A$ if~$\kappa_i=\kappa_j$.
We   briefly review  necessary conditions for diagonal stability of a matrix~$A$
in terms of its principal minors.
%%%
\begin{Proposition} [{\cite[Thm. 2]{barker1978positive}}] \label{prop:ctnec}
Let $A \in \mathbb{R}^{n \times n}$. If there exists a~$D\in \D^{n\times n}$
  such that 
$	 
	D A + A^T D \prec 0,
$	 
	then every principal minor
	of $(-A)$ is positive. 
\end{Proposition}

\updt{Combining this with the Cayley transform  \cite[Thm. 3]{smith1966matrix}  yields the following result.}

\begin{Proposition}[{\updt{see e.g. \cite{kaszkurewicz2012matrix}}}] \label{prop:dtnec}
	Let $A \in \mathbb{R}^{n \times n}$.
	If there exists a~$D\in\D^{n\times n}$ such that 
$	
A^T D A  \prec D,
$
	then  every principal minor  of 
	$
	-(A + I)(A - I)^{-1}
	$
	is positive. 
\end{Proposition}

\updt{
The next three sections describe our main results.}
%%%%%%%%%%%%%%%%%%%%%%%%%%%%%%%%%%
\section{$k$-positivity does not imply diagonal~stability}\label{sec:doesnot}
%%%%%%%%%%%%%%%%%%%%%%%%%%%%%%%%%%%%%%%%%%%%%%%%%%%%%%%%%%%%%%%%%%%%%%%%%%%%%
\updt{Since stable~$1$-positive systems (i.e. positive systems) are diagonally stable, 
a natural question is:
are stable~$k$-positive systems diagonally stable? This section shows that in general the answer is no, both in the DT and CT case.}

Consider the DT  LTI~\eqref{dt_lti} with
	\begin{equation}\label{eq:a7}
	A 
	=\frac{1}{7} \begin{bmatrix}
	-4  & -2  &  1 \\
	1 & -3  & -5  \\
	7  &  1 & -2 \\
	\end{bmatrix}.
    \end{equation}
   It  is straightforward to verify
	that $A$ is Schur, and that~$A^{(2)}$ is $SSR_2$ with $\epsilon_2 = 1$.
	Let $B := -(A + I)(A - I)^{-1}$. Then,
\[	B( \{1,3 \} \vert \{1,3\}) =
	\operatorname{det}\left( \frac{1}{461}
	\begin{bmatrix}
	204  & 140  \\
	497  & 323   \\
	\end{bmatrix} \right) = -\frac{8}{461} <0.
\]
Hence, Proposition~\ref{prop:dtnec} implies that although the DT  LTI
is  stable and strongly 2-positive, it does \emph{not} admit a~DLF.

\begin{Remark}
%%%%%%%%%%%%%%%%%%%%%%%%%%%%%%%
We focus   on DT systems, but  here  we also briefly discuss  the CT case.
The CT LTI~$\dot x=Ax$ 
is called strongly~$k$-positive if its flow maps~$P^k_-\setminus\{0\}$ to~$P^k_+$ 
that is,~$\exp(At) (P^k_-\setminus\{0\})\subseteq P^k_+$ for all~$t>0$. 
	By using Proposition~\ref{prop:ctnec}, we can also prove that $k$-positive CT LTIs are not diagonally stable in general. Consider~$\dot x = A x $ with 
	\begin{equation*}
	A = 
	\begin{bmatrix}
	-21  & 11  &  -14 \\
	18  & -19  &  37  \\
	-49  &  21  & -33 \\
	\end{bmatrix}.
	\end{equation*}
	This system is strongly $2$-positive   (see~\cite{Eyal_k_posi}), 
	and $A$ is Hurwitz. Let $B := -A$. 
	Then
	\[	B( \{2,3 \} \vert \{2,3\}) =
	\operatorname{det}\left( 
	\begin{bmatrix}
	19  & -37  \\
	-21  & 33   \\
	\end{bmatrix} \right) = -150 <0.
	\]	  
	Thus, Proposition~\ref{prop:ctnec} implies that \updt{this system} is \emph{not}~diagonally stable.
\end{Remark}

Summarizing,~stable $k$-positive LTIs are in general not diagonally stable. 
A natural question then is what can be said about the diagonal stability of 
such systems. 

\section{DT $k$-diagonal stability} \label{sec:k-diag}
%%%%%%%%%%%%%%%%%%%%%%%%%%%%%%%%%%%%%%%%%%%%%%
%%%%%%%%%%%%%%%%%%%%%%%%%%%%%%%%%%%%%%%%%%%%%
We begin with defining a new notion called~\emph{$k$-diagonal stability}.
%%%%%%%%%%%%%%%%%%%

\begin{Definition} \label{def:k-diag}
%%%%%%%%%%%%%%%%%%%%%%%%%%%%
Given~$A \in \R^{n\times n}$ and~$k \in [1,n-1]$, let~$r:=\binom{n}{k}$. 
We say that~$A$ is \emph{DT $k$-diagonally stable}
if there exists~$D\in \D^{r \times r}$ such that
\be \label{eq:akds}
 (A^{(k)})^T D A^{(k)}
	\prec D .
\ee
\end{Definition}

Note that Definition \ref{def:k-diag} reduces to   standard DT diagonal stability for~$k=1$, as then~$A^{(1)}=A$ and~$r=\binom{n}{1}=n$. 
The next result is a generalization of Lemma~\ref{le_diag}.
It shows that a $k$-positive DT  LTI is $k$-diagonally stable iff $A^{(k)}$ is Schur.  
%%%%%%%%%%%%%
\begin{Corollary}  \label{cor:k_diag_dt}
%%%%%%%%%%%%%%%%%%%%%%%%%%%%%%%%%%%%55
 Suppose that~$A \in \mathbb{R}^{n \times n}$ is $SR_k$ for some~$k\in[1,n-1]$,
with~$\epsilon_k=1$. Let~$r:=\binom{n}{k}$. Then
 the following statements are equivalent:
	\begin{enumerate}[(a)]
	\item The matrix $A^{(k)}$ is Schur;\label{cond:ksc}
	\item There exists $\xi \in\R^r$ with~$\xi \gg 0$ such that  $A^{(k)}  \xi \ll \xi$;
	\item There exists  $z\in\R^r$ with~$z \gg 0 $ such that    $(A^{(k)})^T z \ll z$;
	\item There exists  $D \in \D^{r \times r}$ such that
	\eqref{eq:akds} holds;
	\item  $(I - A^{(k)})$ is nonsingular and $(I - A^{(k)})^{-1} \geq 0 $. \label{cond:nsin}
	\end{enumerate} 
\end{Corollary} 
%%%%%%%%%%%%%%%%%%%%%%%%%%%%%%%%%%%%%
\begin{comment}
\begin{IEEEproof}
	\item $\prod_{i=1}^k \lambda_i(A) \in [0,1)$.\label{cond:prodeig}
If~$A$ is~$SR_k$ with~$\epsilon_k=1$ then~$A^{(k)}\geq 0$. Thus, 
the equivalence between Conditions~\eqref{cond:ksc}-\eqref{cond:nsin} follows from Lemma~\ref{le_diag}. 

To prove the equivalence with Condition~\eqref{cond:prodeig}, we first pose a stronger assumption, namely, that~$A$ is nonsingular and~$SSR_k$, with~$\epsilon_k=1$.
It then follows from~\cite[Thm.~2]{rola_spect}   
that~$\eta:= \prod_{i=1}^k \lambda_i(A)$ is real and positive. 
If Condition~\eqref{cond:ksc} holds then the modulus of every eigenvalue  of~$A^{(k)}$ is
 smaller than one. Since~$\eta$ is an eigenvalue of~$A^{(k)}$, $\eta \in (0,1)$, so 
Condition~\eqref{cond:prodeig} holds. Conversely, assume that Condition~\eqref{cond:prodeig} holds. Since~$\eta=\rho(A^{(k)})$, we conclude that  Condition~\eqref{cond:ksc} holds. 
This shows that Conditions~\eqref{cond:ksc} and~\eqref{cond:prodeig} are equivalent under the stronger assumption that~$A$ is nonsingular and~$SSR_k$, with~$\epsilon_k=1$. 
Arguing as in~\cite[Section~2.2]{pinkus}, it is not difficult
 to show that given~$A$ that is~$SR_k$, with~$\epsilon_k=1$,
and~$\alpha>0$ there exists a matrix~$B$ that is nonsingular,  $SSR_k$, with~$\epsilon_k=1$,
and~$||A-B||\leq \alpha$. Now a continuity argument  completes the proof of Corollary~\ref{cor:k_diag_dt}.
%%%%%%%%%%%%%
\end{IEEEproof}
\end{comment}

\begin{Remark}\label{rem:frte}
	Note that when these conditions hold we can use the idea described in  Remark~\ref{re_proce} to get an explicit matrix~$D \in \D^{r\times r }$ such that~\eqref{eq:akds} holds. 
\end{Remark}

To demonstrate an   application of Corollary~\ref{cor:k_diag_dt},
we revisit the class of \emph{cyclic DT LTIs}, whose diagonal stability has been analyzed
 in~\cite{wimmer2009diagonal}. 
\begin{Definition}
The matrix 
  $A \in \mathbb{R}^{n \times n}$ is called  cyclic if 
	\begin{equation} \label{eq:cy_dt}
	A = 
	\begin{bmatrix}
	\alpha_1  & \beta_1  &  0 & \cdots &0 \\
	0  & \alpha_2  & \beta_2  & \cdots & 0  \\
	0  &  0  & \alpha_3 & \cdots & 0 \\
	\vdots & \vdots & \vdots & \ddots & \vdots \\
	0 & 0 & 0 & \cdots & \beta_{n-1} \\
	(-1)^{\ell+1} \beta_n & 0 & 0 & \cdots & \alpha_n
	\end{bmatrix},
	\end{equation} 
	with $ \alpha_i, \beta_i \geq 0$, $i = 1, \dots, n$,
	and~$\ell \geq 0 $ is an integer. 
	%%%
	\end{Definition}
We say that the DT LTI~$x(j+1)=Ax(j)$ is cyclic if~$A$ is cyclic. 
Then the dynamics represents a linear chain such that~$x_i(j+1)$ 
depends only on~$x_i(j),x_{i+1}(j)$, and~$x_n(j+1)$ also depends on a feedback connection from~$x_1(j)$. The feedback is negative [positive] if~$\ell$ is even [odd].
%%%%%%%%
 The  next result  shows that such systems are 
DT~$k$-diagonally stable.

\begin{Theorem} \label{thm:cy_srk}
%%%%%%%%%%%%%%%%%%%%%%%%%%%%%%%%%%%%%%%%%
Suppose that~$A$ is cyclic for some~$\ell \in [ 1,  n-1 ]$.
Then~$A$ is $SR_\ell$ with   signature $\epsilon_\ell = 1$.
%%%%%%%%
\updt{Furthermore, if $\ell$ is odd, then~$A$    is 
 DT diagonally stable iff $A$ is Schur.
 If $\ell$ is even,    then $A$ is  DT $\ell$-diagonally stable iff  $A^{(\ell)}$ is Schur.}
%%%%
\end{Theorem} 

\begin{IEEEproof} 
%%%%%%%%%%%%
Pick $\kappa_i, \kappa_j \in Q_{\ell,n}$.
By the Leibniz formula, 
	\begin{align}  \label{eq:akij}
	A(\kappa_i \vert \kappa_j) =& \det(A[\kappa_i \vert \kappa_j])\nonumber \\=& \sum_{\sigma \in \operatorname{pt}(\kappa_j) } \left( \sgn(\sigma) \prod_{s = 1}^\ell  a_{\kappa_{is}, \sigma_s} \right),  
 \end{align}
		where~$\kappa_{is}$ is
		the $s$th element of $\kappa_i$, 
		$ \sigma_s$ is the $s$th element of the permutation $\sigma \in \operatorname{pt}(\kappa_j)$,   $\sgn(\sigma) \in \{-1, 1\}$ denotes the signature of $\sigma$,
	and~$\operatorname{pt}(\kappa_j)$ denotes the set 
	of all~$\ell!$ permutations of the indexes in~$\kappa_j$. 
	For example, if $n=7$, $\ell =3$, and $\kappa_j =\{2,5,7 \}$, then 
	\begin{equation*}\begin{split}
	\operatorname{pt}(\kappa_j) =& \{ \{ 2, 5,7\}, \{ 2, 7,5\}, \{5,2,7\}, \\
	& \{5,7,2\}, \{7,5,2\}, \{7,2,5\} \}.
	\end{split}\end{equation*}
  
	The cyclic structure~\eqref{eq:cy_dt} 
	implies that~$\sgn(\sigma) \prod_{s = 1}^\ell a_{\kappa_{is}, \sigma_s} $ can be
	non-zero
	only in the  following   cases:
	%%%%%%%%%%
	\begin{enumerate} [(i)]
		\item $\kappa_{i \ell} \leq n-1$, and either~$\sigma_s = \kappa_{is}$ or  $\sigma_s = \kappa_{is} + 1$, for all $s \in \{ 1, \dots, \ell \}$;
		\item $\kappa_{i \ell} = n$, $\sigma_\ell = n$, and either~$\sigma_s = \kappa_{is}$ or  $\sigma_s = \kappa_{is} + 1$ for all $s \in \{ 1, \dots, \ell-1 \}$;
		\item $\kappa_{i \ell} = n$, $\sigma_\ell = 1$, and either~$\sigma_s = \kappa_{is}$ or  $\sigma_s = \kappa_{is} + 1$ for all $s \in \{ 1, \dots, \ell-1 \}$.
	\end{enumerate}
	Additionally,  the elements of $\sigma$ are   distinct,
	as~$\sigma \in \operatorname{pt}(\kappa_j)$.
	Since~$\kappa_i$ is an increasing sequence, in
	Cases~$(i)$ and~$(ii)$ the number of inversions
	in~$\sigma$ is zero, so~$\sign(\sigma)=1$. 
	If Case~(iii) holds,
	then~$\sigma_\ell = 1 < \sigma_2 < \dots < \sigma_{\ell-1}$, 
	so~$\sigma$ has~$\ell-1$ inversions. 
	
	Assume that~$\ell$ is even. 
	Then all the entries of~$A$ are non-negative, except perhaps for~$a_{n1}$. 
	 In Cases~(i) or~(ii) we have~$\sgn(\sigma) \prod_{s = 1}^\ell a_{\kappa_{is}, \sigma_s} \geq 0$ since $\sgn(\sigma) =1$, and all the $a_{ij}$s  
	in $\prod_{s = 1}^\ell a_{\kappa_{is}, \sigma_s}$ are non-negative. If Case~(iii) holds, 
	then the  number of inversions in $\sigma$ is $\ell-1$, which is odd, so~$\sgn(\sigma) = -1$. Furthermore,   $a_{n1} \leq 0$ appears in the term~$\prod_{s = 1}^\ell a_{\kappa_{is}, \sigma_s}$.
	Thus, in this case we also have~$\sgn(\sigma) \prod_{s = 1}^\ell  a_{\kappa_{is}, \sigma_s} \geq 0$. Now~\eqref{eq:akij} implies 
	that~$A(\kappa_i \vert \kappa_j) \geq 0$. Since~$\kappa_i, \kappa_j \in Q_{\ell,n}$ 
	are arbitrary, we conclude that~$A$ is $SR_\ell$ with  signature~$\epsilon_\ell =1$.
	The proof for~$\ell$    odd is similar.

Furthermore,
If $\ell$ is even, the results in~\cite{wimmer2009diagonal} show that~$A$ may be Schur yet    not necessarily~diagonally stable. However, since~$A$ is $SR_\ell$,   Corollary~\ref{cor:k_diag_dt} ensures that $A$ is  DT $\ell$-diagonally stable iff  $A^{(\ell)}$ is Schur (which is weaker than the condition~$A$ is Schur).
If $\ell$ is odd, then every entry of~$A$ in \eqref{eq:cy_dt} is non-negative, and Lemma \ref{le_diag} implies that it is 
 DT diagonally stable iff $A$ is Schur.
%%%%%%%%%%%%%%%%%%%%%%%%%%%%%%%%%%%%%%%%%%%%%%%%%%%%%%%%%%%%%
\end{IEEEproof}

\begin{Example}
%%%
			Consider the case~$n=3$, that is,
$
	A = 
	\begin{bmatrix}
	\alpha_1  & \beta_1  &  0  \\
	0  & \alpha_2  & \beta_2    \\
	(-1)^{\ell+1} \beta_3 & 0    & \alpha_3
	\end{bmatrix}.
$ 
%%%
A calculation gives
\[
A^{(2)}= 
	\begin{bmatrix}
	\alpha_1 \alpha_2 & \alpha_1 \beta_2  &  \beta_1\beta_2  \\
	 (-1)^{\ell} \beta_1\beta_3 & \alpha_1\alpha_3  & \alpha_3\beta_1    \\
	(-1)^{\ell} \alpha_2 \beta_3 & (-1)^{\ell} \beta_2\beta_3    & \alpha_2\alpha_3
	\end{bmatrix}.
\]
If~$\ell=1$ then all the entries of~$A$ are non-negative, so~$A$ is~$SR_1$ with signature~$\epsilon_1=1$. 
If~$\ell=2$ then all the entries of~$A^{(2)}$ are non-negative, so~$A$ is~$SR_2$ with signature~$\epsilon_2=1$. 
\end{Example}

The next result shows that DT $k$-diagonal stability, with~$k>1$,
 is a necessary condition for DT diagonal stability. Let~$I_s$   denote the~$s\times s$
identity matrix. 
%%%%%%%%%%%%%%%%%%%%%%%%%%%%%%%%%%%%%%%%%%%%%%%%%%%%%%%%%% 
\begin{Theorem}  \label{thm:k_diag_nec}
If~$A \in \mathbb{R}^{n \times n}$ is DT diagonally stable,
 then $A$ is  DT $k$-diagonally stable for any~$k \in [1, n-1]$.
\end{Theorem}

\begin{IEEEproof}
	Since $A$ is DT diagonally stable, there exists~$P \in \mathbb{D}^{n \times n}$ such that $A^T P A \prec P$. Hence, $P^{-\frac12} A^T P A P^{-\frac12}\prec I_n$,
	so~$P^{-\frac12} A^T P A P^{-\frac12}$ is Schur. 
	Pick~$k\in[1,n-1]$, and let~$r:=\binom{n}{k}$ and~$D := P^{(k)}$. Note that~$D\in \D^{r\times r}$.
	  Lemma~\ref{le_multi} implies  that
	\begin{equation*}
	(P^{-\frac12} A^T P A P^{-\frac12})^{(k)} =
	D^{-\frac12} (A^{(k)})^T D A^{(k)} D^{-\frac12} 
    \end{equation*}
	is also Schur, i.e., $D^{-\frac12} (A^{(k)})^T D A^{(k)} D^{-\frac12} \prec I_r$. We conclude
	that~$(A^{(k)})^T D A^{(k)} 
	\prec D$.
\end{IEEEproof}

Corollary~\ref{cor:k_diag_dt} guarantees
 that a stable $k$-positive DT  LTI  
 is always $k$-diagonally stable with a matrix~$D \in \mathbb{D}^{r \times r}$. 
If there exists~$P \in \mathbb{D}^{n \times n}$ such that $P^{(k)} = D$, then the proof of
Thm.~\ref{thm:k_diag_nec} suggests
 that~$V(z) :=z^T P z$ is a candidate for a~DLF   for the original DT LTI~$x(j+1) = Ax(j)$. However, for any~$k = [2, n-2]$ and~$D \in \mathbb{D}^{r \times r}$, the equation $P^{(k)} = D$ generally does not admit a solution~$P \in \mathbb{D}^{n \times n}$. 
The next result shows that for~$k = n-1$  this equation is always solvable.

\begin{Theorem}\label{thm:diag}
	%%%%%
	For any~$D\in\D^{n\times n}$,     there
	exists a~$P \in\D^{n\times n}$ such that~$P^{(n-1)}=D$.
	%%%%%%
\end{Theorem}

\begin{IEEEproof}
	%%%%%%
	The proof is constructive. 
	The equation~$P^{(n-1)}=D$~can be written as
	\be \label{eq:sqn-1}
	\prod_{s \in \kappa_q}p_s = d_q  ,\quad  q = 1, \dots, n
	\ee
	where $\kappa_1, \dots, \kappa_n \in Q_{n-1,n}$, and $p_i, d_i$ denote the $i$th diagonal entry of $P$ and $D$, respectively.
	For any~$s\in[1,n]$, let~$ j(s)  $ be the single element
	in the set of indexes~$[1,n] \setminus \kappa_s  $.
	A lengthy but straightforward computation shows that  
	the solution of~(\ref{eq:sqn-1})   is 
	\be\label{eq:gtpos}
	p_s = \frac{\prod_{q \in \kappa_s} d_q^{\frac{1}{n-1}}}{d_{j(s)}^{\frac{n-2}{n-1}}}.
	\ee
	Since~$d_i>0$ for any~$i$,  this implies that~$p_s >0$ for any~$s$.
	%%%%%%
\end{IEEEproof}

The following example shows that how the above results
 can be utilized to construct a DLF for an $(n-1)$-positive DT LTI.

\begin{Example}
	%%%%%%%%%%%%%%%%%%%%%%%%%%%%%%%%%%%
	Consider the DT LTI $x(j+1) = Ax(j)$ with
	\begin{equation}
	A = \frac{1}{8}
	\begin{bmatrix}
	-4  & -2  &  0 \\
	0  & -3  & -5  \\
	7  &  0  & -2 \\
	\end{bmatrix} .
	\end{equation}
	A calculation shows that $A$ is Schur.
	Since the entries of~$A$ have different signs, 
	we cannot use  Lemma~\ref{le_diag} to conclude that~$A$
	admits a~DLF. However,~$A$
	is~$SSR_2$ with~$\epsilon_2 = 1$. Hence, Corollary~\ref{cor:k_diag_dt} implies that there exists~$D \in \D^{3 \times 3}$ such that $  (A^{(2)})^T D  A^{(2)}   \prec D$. According to Remark~\ref{rem:frte}, 
	one such $D$ can be obtained as 
	$
	D = \diag \left(\frac{23}{21}, \frac{13}{8},\frac{7}{13} \right).
	$
	Using Theorem~\ref{thm:diag} to solve $P^{(2)} = D$ gives
	$P = \diag \left
	(\sqrt{\frac{3887}{1176}}, \sqrt{\frac{184}{507}}, \sqrt{\frac{147}{184}} \right).
	$
	It is straightforward to verify that $A^T P A \prec P$. Thus, we were able
	to build a DLF
	for~$A$.
\end{Example}

\section{Applications to Nonlinear Dynamical Systems} \label{sec:non}
As mentioned in the introduction,   DT diagonal stability of~$A$ implies that certain nonlinear DT systems are also stable.
 A natural question is what are the implications of DT $k$-diagonal stability for nonlinear systems?
In this section, we describe a new  class
  of DT nonlinear system whose dynamics can be analyzed by 
	exploiting $k$-positivity and wedge products.
We first define a special kind of nonlinear mappings.
%%%%% 
 \begin{Definition} \label{def:kmap}
 	Let $\S\subseteq\R$ with~$0\in \operatorname{int}\S$.
 	Define $\phi:\S^n \to \R^n$ by $\phi(x) := \begin{bmatrix} \phi_1(x_1)& \dots &\phi_n(x_n) \end{bmatrix}^T$, where every~$\phi_i:\S\to\R$
	is a continuous scalar function such that~$\phi_i(s)=0$ holds only for~$s=0$.
	Pick~$k\in[1,n-1]$ and let~$r:=\binom{n}{k}$.
	We say that~$\phi$ is \emph{$k$-content preserving} if for any~$a^1,\dots,a^k \in \S^n$ we have that 
	\begin{equation}  \label{eq:kvolpre}
			%%%
			\begin{cases}
	  q_i = 0,  & \text{ if } p_i = 0, \\
	    | q_i |  \in (0,\vert p_i \vert]  & \text{ if } p_i \neq 0,								
		%%%%
	\end{cases}\end{equation}
	for all $i=1,\dots,r$, where~$q:=\wedge_{j=1}^k \phi(a^j)$ and~$p:= \wedge_{j=1}^n a^j$.
%%%%%%%%%%%%
\end{Definition}

\begin{Example}
For~$k=1$ we have~$p=a $, $q=\phi(a )$, so~\eqref{eq:kvolpre} reduces to $\phi_i(0) = 0$, and $ 0<|\phi_i(a_i)|\leq |a_i|$ for $s \neq 0$.
%%%
   For~$k=2$, pick~$a,b \in \R^n$, and
let~$p:=a\wedge  b$, $q:=\phi(a) \wedge \phi(b)$.  
Then   $p_1=  a_1 b_2-a_2 b_1   $ and~$q_1=\phi_1(a_1) \phi_2(b_2)-\phi_2(a_2) \phi_1(b_1 ) $.
Thus, for~$i=1$, \eqref{eq:kvolpre} yields
\be\label{eq:jdhe1}
 \vert	\phi_1(a_1) \phi_2(b_2)-\phi_2(a_2) \phi_1(b_1 ) \vert \leq \vert a_1 b_2-a_2 b_1 \vert.							%%%
\ee
%%%%%%%%%%%%%%%%%%%%%%%%%%%%%%%%%%%%%%%
(The equations for other values of~$i$ are similar.) 
For~$a_2=0$ this gives~$(\phi_1 (a_1) \phi_2 (b_2 ))^2 \leq   a_1^2 b_2 ^2$.
If~$a= \alpha b$, with~$\alpha \in \R\setminus\{0\}$, 
then~\eqref{eq:jdhe1} becomes
\[
							(\phi_1(\alpha b_1) \phi_2(b_2)-\phi_2(\alpha b_2) \phi_1(b_1 ))^2\leq 0,
\]
that is,
\[
							 \phi_1(\alpha b_1) \phi_2(b_2) = \phi_2(\alpha b_2) \phi_1(b_1 ),
\]
and for~$b_2\not =0$ this becomes the homogeneity condition
\[
							\frac{ \phi_1(\alpha b_1) }{\phi_2(\alpha b_2) }  =  \frac{ \phi_1(b_1 )}{ \phi_2(b_2)}.
\]
As a specific example, take~$\S=[-1/2,1/2]$ and~$\phi_i(s)=s^2$, for all~$i$. 
Then it is not difficult to show that~\eqref{eq:jdhe1}
 holds for any~$a_i,b_i\in\S$, so this function is~$2$-content preserving on $\S$.  
%%%%%%%%%%%%%%%%%%%%%%%%%%%%%%%%%5
\end{Example}

We can now state the main result in this section.
\begin{Theorem}   \label{thm:k_diag_dt_non}
%%%%%%%%%%%%%%%%%%%%%%%%%%%%%%%%%%%%
Suppose that~$A\in\R^{n\times n}$ is DT~$k$-diagonally stable for some~$k\in [1,n-1]$. 
Consider the~DT nonlinear system
\be\label{eq:dtn}
x(j+1)=A\phi(x(j)),
\ee
 where~$\phi(x)=\begin{bmatrix}
\phi_1(x_1)&\dots&\phi_n(x_n)\end{bmatrix}^T$ is~$k$-content preserving on the state-space~$\S^n$ of~\eqref{eq:dtn}. 
%%%%%%%%%%%%%%%%%%%%%%%%%%%%%%%%%%%%%%%%%%%%%%%%%%%%%%%%%%%%%%%%%%%%% 
 For~$a^1,\dots,a^k\in\S^n$, let
\be \label{eq:yde}
 y(j) = y(j;a^1,\dots,a^k):=  \wedge_{i=1}^k x(j,a^i), 
\ee
where~$x(j,a)$ is the solution of~\eqref{eq:dtn} at time~$j$ with~$x(0)=a$. Then
\be \label{eq:ydin}
	y(j+1) =A^{(k)} \wedge_{i=1}^k \phi(x(j,a^i)),
\ee
and this nonlinear dynamical system is diagonally stable.
%%%%%%%%%%%%%%%%%%%%%%%%%%%%%%%%%%%%%
\end{Theorem} 

\begin{IEEEproof}
%%%%%%
	By   \eqref{eq:yde},
\begin{align*}
%%%%%
						y(j+1)&=  \wedge_{i=1}^k x(j+1,a^i)\\
		      &=A\phi(x(j,a^1))\wedge \dots\wedge A\phi(x(j,a^k))\\
									&= \begin{bmatrix} A\phi(x(j,a^1))&  \dots& A\phi(x(j,a^k)) \end{bmatrix}^{(k)}\\
								%%	&= (A \begin{bmatrix}  \phi(x(j,a^1))& \dots & \phi(x(j,a^k)) \end{bmatrix})^{(k)}\\
				&= A^{(k)} \begin{bmatrix}  \phi(x(j,a^1))& \dots&  \phi(x(j,a^k)) \end{bmatrix}^{(k)},
			%%%%%
\end{align*}
and this proves~\eqref{eq:ydin}. 
Since~$A$ is DT~$k$-diagonally stable, there exists~$D\in\D^{r\times r}$ such
 that~\eqref{eq:akds} holds. 
Define~$V:\R^r\to\R_+$ by~$V(z):=z^T D z $, and let~$\triangle V(j):=V(	y(j+1) ) -V(y(j) )$.
 Then
		\begin{align}\label{eq:vdify}
%%%%%
		\triangle V(j) 
		=& 
		( \wedge_{i=1}^k \phi(x(j,a^i)) )^T (A^{(k)})^T D A^{(k)} \wedge_{i=1}^k \phi(x(j,a^i)) \nonumber \\
		&- ( \wedge_{i=1}^k x(j,a^i) )^T   D   \wedge_{i=1}^k x(j,a^i)  .
%%%%%
\end{align}
	 	%%%%%
	  Definition~\ref{def:kmap}
		implies that:
		\begin{align*} 
	 (\wedge_{i=1}^k \phi(x(j,a^i)) )^T & D  \wedge_{i=1}^k \phi(x(j,a^i)) \\&
	  \leq (  \wedge_{i=1}^k x(j,a^i)  )^T   D   \wedge_{i=1}^k x(j,a^i) ,
		\end{align*}
		and combining this with~\eqref{eq:vdify} gives
				$\triangle V(j) \leq   (\wedge_{i=1}^k \phi(x(j,a^i)) )^T ( (A^{(k)})^T D A^{(k)}- D)  \wedge_{i=1}^k \phi(x(j,a^i)).$
	We conclude that~$V(y(j+1) ) -V(y(j) )\leq 0$, with equality only when~$y(j)=0$.
%%%%%%%%%%%%
\end{IEEEproof}

%%%%%%%%%%%%%%%%%
Note   that the existence of a~$D \in \D^{r \times r}$ that satisfies~\eqref{eq:akds}
plays a  crucial role in the proof.  
%%%%%%%%%%%% 

Theorem~\ref{thm:k_diag_dt_non} implies that the $k$-content of the parallelotope induced by $x(j, a^i)$, $i =1, \dots, k$, converges to zero asymptotically. For~$k=2$
 this means that any two  trajectories of~\eqref{eq:dtn}   converge to a line, i.e., to a one-dimensional subspace. In particular,
this ensures that the dynamics of \eqref{eq:dtn} has no nontrivial limit cycles. 

\begin{comment}
\begin{Remark}
By Corollary~\ref{cor:k_diag_dt}, Theorem~\ref{thm:k_diag_dt_non} holds if $A$ is $SR_k$ and $A^{(k)}$ is Schur. However, it is worthy of pointing out that $A\frac{\partial }{\partial x} \phi(x)$ is not necessarily $SR_k$ even if $A$ is $SR_k$. Recall the definition of cooperative systems \cite{hlsmith}, the notion of $k$-cooperative systems can be introduced. That is, consider the nonlinear system
\begin{equation} \label{eq:xfnon}
x(j+1) = f(x(j)),
\end{equation}
where $x \in \R^n$, and $f: \R^n \to \R^n$ is a continuous and differentiable mapping. We say that \eqref{eq:xfnon} is $k$-cooperative if $\frac{\partial }{\partial x}f(x)$ is $SR_k$ for all $x \in \R^n$. Clearly, the nonlinear systems \eqref{eq:dtn} are not $k$-cooperative in general.
\end{Remark}
\end{comment}

\begin{Corollary}   \label{prop:cyc_non}
	%%%%%%%%%%%%%%%%%%%%%%%%%%%%%%%%%%%%
	Consider the~DT nonlinear system:
	\be\label{eq:cyc_non}
	x(j+1)=A\phi(x(j)),
	\ee
	where $A$ is cyclic 
	for some~$\ell \in [1,n-1]$, and~$\phi(x)=\begin{bmatrix}
	\phi_1(x_1)&\dots&\phi_n(x_n)\end{bmatrix}^T$ is~$\ell$-content preserving on the state-space~$\S^n$
	of~\eqref{eq:dtn}. 
%%%%%%%%%%%%%%%%%%%%%%%%%%%%%%%%%%%%%%%%%%%%%%%%%%%%%%%%%%%%%%%%%%%%% 
	For~any $a^1,\dots,a^\ell \in\S^n$, let
	$
	y(j) = y(j;a^1,\dots,a^\ell):= \wedge_{i=1}^\ell  x(j,a^i).
	$ 
	Then 
	\be \label{eq:cyc_ydin}
	y(j+1) =A^{(\ell)} \wedge_{i=1}^\ell \phi(x(j,a^i)),
	\ee
	and if $ \vert \prod_{i=1}^\ell
		\lambda_i(A) \vert < 1$ then~\eqref{eq:cyc_ydin}
	 is diagonally stable.
	%%%%%%%%%%%%%%%%%%%%%%%%%%%%%%%%%%%%%
\end{Corollary} 

\begin{IEEEproof}
%%%%%%%%%%%%%%%%%%%%%%%%%%%%%%%%%%%%%%%%%%%%
	  By Theorem~\ref{thm:cy_srk},
		 $A$   is $SR_\ell$ with $\epsilon_\ell =1$.  By Corollary~\ref{cor:k_diag_dt}, 
		$A$ is DT diagonally stable iff 
		 $A^{(\ell)}$ is Schur, that is, iff~$ \vert \prod_{i=1}^\ell \lambda_i(A) \vert < 1$. 
		Applying Theorem~\ref{thm:k_diag_dt_non} completes the proof.
%%%%%%%%%%%%
\end{IEEEproof}

\begin{figure}[t]
	\begin{center}
		%%%%%%%%%%%%%%%%%%%%%%%[width=8cm,height=12cm] scale=1.0
		\includegraphics[scale=0.4]{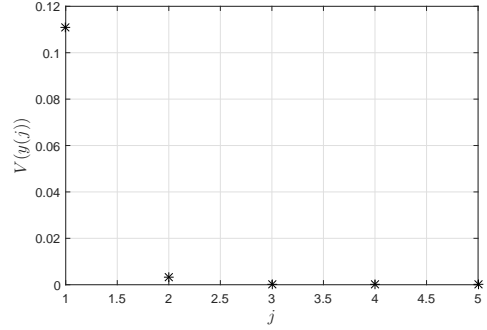}
		\caption{$V(y(j))$ as a function~$j \in [1,5]$  in Example~\ref{exa:cyc}. }
		\label{fig:gap}
		%%%%%%%%%%%%%%%%%%%%%%%%%%
	\end{center}
\end{figure}

\begin{Example}\label{exa:cyc}
	%%%%%%%%%%%%%%%%%%%%%%%%%%%
	Consider the DT nonlinear system~\eqref{eq:dtn}
	with~$n=3$, $A=\begin{bmatrix}  0.1 &1.9& 0  \\
	0& 0.05 & 1.95  \\ -0.01 & 0& 2.01 \end{bmatrix} $, and~$\phi_i(s)=s^2$, $i = 1,2, 3$, that is,~$\phi(x)=\begin{bmatrix} x_1^2 & x_2^2 & x_3^2\end{bmatrix}^T$. 
	Let~$\S:=[-1/2,1/2]$. It is not difficult to show that~$\S^3$ is an invariant set of the dynamics. For example,
	$x_1(j+1)= 0.1 x_1^2(j) + 1.9 x_2^2 (j)    $.
	If~$x(j) \in \S^3$ then~$ x_i^2(j) \in [0,1/4]$ and this implies that~$x_1(j+1)\in \S$. 
	The matrix~$A$ is \emph{not} Schur, as~$\rho(A)=2$.
	However,~$A^{(2)}$ is Schur, and also~$ A^{(2)} \geq 0$, i.e.~$A$ is~$SR_2$ with~$\epsilon_2=1$.
	We use the idea  described in
	Remark~\ref{re_proce}  to get a~$D$ such that~\eqref{eq:akds} holds. 
	Here~$n=3$ and~$k=2$, so~$r=\binom{n}{k}=3$. 
	Denote~$1_3:=\begin{bmatrix} 1&1&1\end{bmatrix}^T$, and let 
	\begin{align*} 
	\xi& := (I - A^{(2)})^{-1} 1_3 ,\\
	z&:=(I-(A^{(2)})^T)^{-1} 1_3 , \\
	D&:=\diag( z_1/ \xi_1,  z_2/ \xi_2,z_3/\xi_3) .
	%%%%%%%%%%%%%%%%%%%%%%%%%%%%%
	\end{align*}
	Fig.~\ref{fig:gap} depicts~$V(y(j))=y^T(j)Dy(j)$, as a function of~$j$,
	where~$y(j):=x(j,a^1)\wedge x(j,a^2)$, for the initial conditions
	$a^1=(1/2)1_3$,
	$a^2=\begin{bmatrix}  -1/2& 1/2& 0.4 \end{bmatrix}^T $. Note that~$a^1,a^2\in \S^3$. 
	As expected,~$V(y(j))$ decreases with~$j$. 
	
	If we take $a:=(1/2)1_3  $,   $b \in \S^3$, then
	\begin{align*}
	y(j)&=x(j,a)\wedge x(j,b)\\
			&=a \wedge x(j,b)\\
			&=(1/2)\begin{bmatrix} x_2-x_1& x_3-x_1& x_3-x_2  
			\end{bmatrix}^T,
	\end{align*}
	where~$x_i:=x_i(j,b) $.
	Thus,~$4V(y(j))$ is equal to
	$$   d_1 (x_2 -x_1)^2 +   d_2 (x_3 -x_1)^2 +   d_3 (x_3 -x_2)^2  ,$$
	where~$d_i$ is the $i$th diagonal entry of~$D$.
	Since we already know that this
	function  converges   
	to zero, every trajectory converges
	to the line spanned by~$1_3$.
	%%%%%%%%%%%%%%%%%%
\end{Example}

\section{Conclusion} \label{sec:con}
%%%%
Diagonal  stability is an important property of positive LTIs. 
$k$-positive LTIs are a generalization of positive LTIs and so a natural question is whether stable $k$-positive LTIs are also diagonally stable. We showed that in general the answer is no. 

We then defined the new notion of DT~\emph{$k$-diagonal stability} and showed how it can be used to generalize the idea that diagonal stability of an LTI implies the  stability of a certain nonlinear dynamical system. These results admit a clear geometric interpretation using  the wedge product. We demonstrated our results for a class of nonlinear systems that include a cyclic matrix in their dynamics. 
 
Due to space limitations, we focused here on DT systems. The CT case       may be an interesting topic for further research. 

\begin{comment}
\updt{\subsection*{Acknowledgements}
We are grateful   to the AE and the anonymous  reviewers for many helpful comments.
}
\end{comment}

%%%%%%%%%%%%%%%%%%%%%%%%%%%%%%%%%%%%%%%%%%%%%%%%%%%%%%%%%%%%

\begin{comment}
***** MATLAB PROGRAM **********************
***** DO NOT DELETE *********************

%%%FILE DIAG_STAB.m
clear all
hold off; 
A=1.7* [.2, 1, -.1 ;  .1, 1, 0  ; 0 , 1, 0.1]; 
A=[.1 1.9 0; 0 1/16 31/16; 0 0 2];
perms=nchoosek(1:3, 2) ; 
B=zeros(3,3);
for i=1:3
    for j=1:3
        alp=perms(i,:);beta=perms(j,:);
    B(i,j)=   det( A(  alp,beta  )   ) ;     %%% B is A^{(2)} 
    end
end
zeta=inv( eye(3)-B) *[1;1;1];
z=inv(eye(3)-B') *[1;1;1];
D=[z(1)/zeta(1) 0 0 ; 0 z(2)/zeta(2) 0 ; 0 0 z(3)/zeta(3)];
a=[1/2 1/2 1/2]'; 
b=[-1/2  1/2 .4]'; 
sol=a;
yol=b;  %%% these are the two solutions each emnating from a different initial condition
for i=1:4
    y=[sol(1)*yol(2)-sol(2)*yol(1) ......
       sol(1)*yol(3)-sol(3)*yol(1) ......
       sol(2)*yol(3)-sol(3)*yol(2)]' ;
   v=y'*D*y;
   r=plot(i,v,'*k','markersize',8);hold on;
   sol=A*(sol.^2)
   yol=A*(yol.^2);
end
       grid on;xlabel('$j$','interpreter', 'Latex', 'fontsize', 16);

%%%%%%%%%%%%%%%%%%%%%%%%%%%%%%%%%%%%%%%%%%%%%%
\end{comment}

%%%%%%%%%%%%%%%%%%%%%%%%%%%%%%%%%%%%%%%%%%%%%%%%%%%%%%%%%%%%
\bibliographystyle{IEEEtran}
\bibliography{refs}
 \end{document}